\documentclass[12pt]{article}
\begin{document}
\input epsf
\setlength{\baselineskip}{1.6\baselineskip}
\newtheorem{lemma}{Lemma}
\newtheorem{cor}{Corollary}
\newtheorem{deff}{Definition}
\newtheorem{pro}{Proposition}
\newtheorem{theorem}{Theorem}
\newtheorem{fig}{Figure}
\newcounter{reno}
\renewcommand{\theequation}{\thesection.\arabic{equation}}
\newcommand{\all}{\renewcommand{\include}{\input}}
\newcommand{\bfb}[1]{
            \mbox{\boldmath $ #1 $}}

\title{An Algorithm to Estimate Monotone Normal Means and its
Application to Identify the Minimum Effective Dose}
\author{By WEIZHEN WANG\\Department of Mathematics and
Statistics\\Wright State University\\Dayton, OH 45435, USA
\\email:wwang$@$math.wright.edu\\
JIANAN PENG\\Department of Mathematics and Statistics\\ Acadia
University,\\ Wolfville, NS  B4P 2R6\\Canada\\ E-MAIL:
jianan.peng@acadiau.ca}
\date{December 29, 2007}
\maketitle
\thispagestyle{empty}

\begin{abstract}
\setlength{\baselineskip}{1.4\baselineskip}

In the standard setting of one-way ANOVA with normal errors, a new
algorithm, called the Step Down Maximum Mean Selection Algorithm
(SDMMSA), is proposed to estimate the treatment means under an
assumption that the treatment mean is nondecreasing in the factor
level. We prove that i) the SDMMSA and the Pooled Adjacent Violator
Algorithm (PAVA), a widely used algorithm in many problems, generate
the same estimators for normal means, ii) the estimators are the
mle's, and iii) the distribution of each of the estimators is
stochastically nondecreasing in each of the treatment means. As an
application of this stochastic ordering, a sequence of null
hypotheses to identify the minimum effective dose (MED) is
formulated under the assumption of monotone treatment(dose) means. A
step-up testing procedure, which controls the experimentwise error
rate in the strong sense, is constructed. When the MED=1, the
proposed test is uniformly more powerful than Hsu and Berger's
(1999).

\bigskip

\noindent {\it Some key words}: Closed test method; Experimentwise
error rate; Maximum likelihood estimator; Step-up tests.

\end{abstract}

\section{Introduction.}\label{intro}
A situation frequently encountered in  dose-response studies is
identifying the minimum effective dose ($MED$). The $MED$ is defined
as the lowest dose such that the mean response is better than that
of a zero-dose control by a clinically significant difference.
Finding the MED is important since high doses often turn out to have
undesirable side effects.

Consider the one-way layout model
\begin{equation}
\label{model} Y_{ij}=\mu_i+\varepsilon_{ij}
\end{equation}
 for $i=0,...,k, j=1,...,n_i$, where $\mu_i$'s are the unknown response
means at different dose levels and $\varepsilon_{ij}\sim
N(0,\sigma^2)$ are the independent errors with an unknown variance.
The parameter space is
\begin{equation}
\label{para} H=\{\underline{\mu}=(\mu_0,\mu_1,...,\mu_k): \mu_1 \leq
...\leq \mu_k\}
\end{equation}
(here, for simplicity, we omit $\sigma$ in $H$), and the sufficient
statistics are the sample means, $\bar{Y}_{i}$, and the mean squared
error, denoted by $S^2$. Assume $i=0$ is the control group. One goal
is to find the smallest positive integer $N$ satisfying $\mu_{N}>
\mu_0+\delta$ for a clinically significant difference constant
$\delta \geq 0$. We call $N$ the minimum effective dose (MED).
Determination of the MED usually is done by step-down test
procedures, see Williams (1971), Ruberg (1989), Tamhane, Hochberg,
and Dunnett (1996),  Hsu and Berger (1999), and Hellmich and
Lehmacher (2005), among others. Tamhane, Hochberg and Dunnett (1996)
indeed proposed a step-up procedure SU1P to identify the MED. The
SU1P procedure is based on the step-up procedure of Dunnett and
Tamhane (1992), which controls the experimentwise error rate only
for balanced designs. However, Dunnett and Tamhane (1995)'s step-up
procedure for unbalanced designs case cannot control the
experimentwise error rate. Liu (1997) proposed a method of
calculating the critical values of the step-up procedure by Dunnett
and Tamhane (1995). The SU1P procedure does not make use of the
monotonicity, therefore its power should not be high. Intuitively it
seems that step-down procedures infer a larger dose as the MED.
Therefore, it is of interest to have a step-up procedure to use the
monotonicity to increase its power.

To derive a test of level-$\alpha$, one needs to find an appropriate
statistic and its least favorable distribution in the null
hypothesis. Thus a stochastic ordering for the test statistic is
needed. The desired statistic, the estimator of $\mu_i$, should be:
a) nondecreasing in $i$, b) and is also nondecreasing in each of
$\bar{Y}_j$'s. The PAVA algorithm generates the estimators that
achieve a). However, it is difficult to show b) directly for these
estimators using the PAVA. The PAVA was first proposed by Ayel,
Brunk, Ewing, Reid and Silverman (1955), and was introduced to
estimate the monotone proportions in independent binomial
experiments. Surprisingly, it has many applications in normal,
Poisson and multinomial distributions, etc. See more details in
Robertson, Wright and Dykstra (1988). The PAVA is an iterative
algorithm, each step is very simple to implement, however, it does
not have a closed form for the final estimator. Hence, it is hard to
establish analytic properties for the estimator. Notice these, a new
algorithm, the SDMMSA, is proposed to overcome the drawbacks. We
will show that the two algorithms yield the same estimators and each
estimator is a monotone function of each $\bar{Y}_i$. The second
fact is critical to determine the least favorable distribution in
the null hypothesis space.

The rest of the article is organized as follows. Section 2 provides
a new algorithm to generate estimators for $\mu_i$'s and discusses
their analytic properties. In particular,  a stochastic ordering for
the distributions of the proposed estimators is established. In
Section 3, one application of the stochastic ordering established in
Section 2 is given to identify the MED. A step-up multiple test
procedure that controls the experimentwise error rate in the strong
sense is provided by constructing a sequence of increasing rejection
regions of level-$\alpha$ for each null hypothesis in
$(\ref{null})$, and the proposed procedure is illustrated on a real
data set. Section 4 concludes with some discussion.

\section{A new algorithm to construct the estimator of $\mu_i$ and some analytic results.}

In this section, an estimator of $\mu_i$, denoted by $\hat{\mu}_i$,
for any integer $i\in [1,k]$ under $H$ is first constructed
iteratively. Then three facts are established: $\hat{\mu}_i$ is the
same as the estimator generated by the PAVA; $\hat{\mu}_i$ is the
mle under $H$, and the distribution of $\hat{\mu}_i$ is
stochastically non-decreasing in each $\mu_j$.

\subsection{A new iterative algorithm to construct  $\hat{\mu}_i$.}

Let
\begin{equation}
\label{combinedmean}  n_{i,j}=\sum_{h=i}^j n_h, \,\
\bar{Y}_{i,j}=\frac{\sum_{h=i}^j n_h\bar{Y}_h}{ n_{i,j}}, \,\
\forall 1\leq i\leq j\leq k
\end{equation}
be the sample size and the sample mean of a combined sample  of
treatments $i$ through $j$, respectively.

Step 1). We construct $\hat{\mu}_i$ starting from $i=k$ using the
data set $\{(\bar{Y}_i, n_i)\}_{i=1}^k$. Let
\begin{equation}
\label{step1set} A_1=\{j: \bar{Y}_{j,k}=max_{\{1\leq j'\leq
k\}}\{\bar{Y}_{j',k} \}   \}
\end{equation}
be a subset of $\{1,...,k\}$($A_1$ contains a single element with
probability one), and let
\begin{equation}
\label{step1index} i_1= min\{A_1\}.
\end{equation}
Then
\begin{equation}
\label{step1mle} \hat{\mu}_{i}\stackrel{def}{=}\bar{Y}_{i_1,k},\,\
\forall \,\ i\in [i_1,k].
\end{equation}
If $i_1=1$, then all $\hat{\mu}_{i}$'s are defined and stop;
otherwise go to the next step.
 Step 2). Note in this step $ i_1-1
\leq k-1$. Repeat Step 1 but using the data set
$\{(\bar{Y}_i,n_i)\}_{i=1}^{i_1-1}$. i.e., let
\begin{equation}
\label{step2set} A_2=\{j: \bar{Y}_{j,i_1-1}=max_{\{1\leq j'\leq
i_1-1\}}\{\bar{Y}_{j',i_1-1} \}   \}
\end{equation}
be a subset of $\{1,...,i_1-1\}$, and let
\begin{equation}
\label{step2index} i_2= min\{A_2\}.
\end{equation}
Then
\begin{equation}
\label{step2mle}
\hat{\mu}_{i}\stackrel{def}{=}\bar{Y}_{i_2,i_1-1},\,\ \forall \,\
i\in [i_2,i_1-1].
\end{equation}
If $i_2=1$, then all $\hat{\mu}_{i}$'s are defined and stop;
otherwise repeat this process for a number of times, say $h$ times,
until $i_h=1$. Such an integer $h$ exists, because $i_j$ strictly
decreases in $j$. Since $i_h=1$, then all $\hat{\mu}_{i}$'s are
defined and the construction on $\hat{\mu}_i$'s is complete. We name
this the step-down-maximum-mean-selection algorithm (SDMMSA).

\noindent{\bf Remark 1}.  There exists partition, $\cup_{u=1}^h
[i_{u},i_{u-1}-1]$, for $\{1,...,k\}$ with
$i_0-1\stackrel{def}{=}k$. Following the construction of
$\hat{\mu}_i$, each
 $\hat{\mu}_i$ is the sample mean of a combined sample of treatment(s) belonging to the interval in
 the partition that includes treatment $i$.
 Also
 $\hat{\mu}_i$ is constant in $i$ on each integer interval $[i_{u}, i_{u-1}-1]$ for $u=1,...h$, as shown in
  ($\ref{step1mle}$) and ($\ref{step2mle}$),
 and $\hat{\mu}_i$  is strictly increasing when $i$ moves from $[i_{u},
 i_{u-1}-1]$ to $[i_{u'}, i_{u'-1}-1]$ for $u>u'$, as shown in ($\ref{step1set}$) and ($\ref{step2set}$). Therefore,
 $\hat{\mu}_i$ is nondecreasing in $i$ for $i\in [1,k]$.
 ~\raisebox{.5ex}{\fbox{}}

\begin{lemma}
\label{mon} For partition $\cup_{u=1}^h [i_{u},i_{u-1}-1]$ given in
Remark 1, $\bar{Y}_{i_u-1}<\bar{Y}_{i_u}$ for any $u\in [1,h]$.
\end{lemma}

\noindent{\bf Proof}. Since the SDMMSA repeats itself in each step,
without loss of generality, we only need to prove Lemma~\ref{mon}
for $u=1$. i.e., $\bar{Y}_{i_1-1}<\bar{Y}_{i_1}$.

Suppose this is not true, i.e.,
 $\bar{Y}_{i_1} \leq \bar{Y}_{i_1-1}$.  Note
$\bar{Y}_{i_1,k} \geq \bar{Y}_{i_1+1,k}$ by the definition of $i_1$,
then
\begin{equation}
\label{contradiciton2}
 \bar{Y}_{i_1} \geq
\bar{Y}_{i_1+1,k}.
\end{equation} Similarly,  $\bar{Y}_{i_1-1} <
\bar{Y}_{i_1,k}$ is true due to $\bar{Y}_{i_1-1,k} <
\bar{Y}_{i_1,k}$.  Therefore, $\bar{Y}_{i_1} \leq \bar{Y}_{i_1-1} <
\bar{Y}_{i_1,k},$ which implies
$$ \bar{Y}_{i_1}<\bar{Y}_{i_1+1,k},$$ a contradiction to ($\ref{contradiciton2}$).
 ~\raisebox{.5ex}{\fbox{}}
\bigskip

\noindent{\bf Remark 2}. If $\bar{Y}_i$ is nondecreasing in
$i\in[1,k]$, then the partition, $\cup_{u=1}^h [i_{u},i_{u-1}-1]$,
for $\{1,...,k\}$ given in Remark 1 satisfies i) $\bar{Y}_i$  is
constant when $i\in[i_{u},
 i_{u-1}-1]$, and ii) is strictly increasing when $i$ moves from $[i_{u},
 i_{u-1}-1]$ to $[i_{u'}, i_{u'-1}-1]$ for $u>u'$. Therefore,
 $\hat{\mu}_i=\bar{Y}_i$  for $i\in [1,k]$.
 ~\raisebox{.5ex}{\fbox{}}

\bigskip

\noindent{\bf Example 1}. Consider the data in Table 1, taken from
Ruberg (1995). There are nine $(k=9)$ active dose groups and a zero
dose control group with six $(n_i=6, i=0,...,9)$ animals/group in
the experiment. Following Step 1, we obtain $i_1=9$ and then
$\hat{\mu}_9=\bar{Y}_9$; following Step 2, we obtain $i_2=6$, then
$\hat{\mu}_8=\hat{\mu}_7=\hat{\mu}_6$ and is equal to
$\bar{Y}_{6,8}=73.77$, the sample mean of the combined sample for
$i=8,7,6$. The construction of all $\hat{\mu}_i$'s ends at Step
7(=h) and their values are reported in Table 1. The partition given
in Remark 1 is now
$$[9]\cup [6,8] \cup [5]\cup [4]\cup [3]\cup [2]\cup [1],$$
with a notation of $[i]=[i,i]$.
 ~\raisebox{.5ex}{\fbox{}}

\bigskip
\subsection{The relationship between $\hat{\mu}_i$,   $\hat{\mu}_i^{mle}$ and $\hat{\mu}_i^{pava}$.}

So far, the estimator of $\mu_i$ under $H$ typically is obtained
following the pooled-adjacent-violators algorithm(PAVA, described
later), for example, see Barlow, Bartholomew, Bremner and Brunk
(1972), Robertson, Wright, and Dykstra (1988), and Silvapulle and
Sen (2005). Denote this estimator by $\hat{\mu}_i^{pava}$. Now we
show that $\hat{\mu}_i^{mle}=\hat{\mu}_i=\hat{\mu}_i^{pava}$ in
Theorem~\ref{sdmmsamle} and Theorem~\ref{mle} below.

\begin{theorem}
\label{sdmmsamle} Let $\hat{\mu}_i^{mle}$ be the maximum likelihood
estimator for $\mu_i$ under $H$. Then
$\hat{\mu}_i=\hat{\mu}_i^{mle}, \forall i \in [1,k]$. Therefore,
$\hat{\mu}_i=\hat{\mu}_i^{mle}, \forall i \in [1,k]$.
\end{theorem}

\noindent{\bf Proof of Theorem~\ref{sdmmsamle}}. Taking the log
transformation on the joint pdf of $Y_{ij}$, it is easy to see that
$\hat{\mu}_i^{mle}$ minimizes
$$f(\mu_1,...,\mu_k)=\sum_{i=1}^k n_i (\bar{Y}_i-\mu_i)^2=\sum_{u=1}^h
[\sum_{j=i_u}^{i_{u-1}-1} n_j
(\bar{Y}_j-\mu_j)^2]\stackrel{def}{=}\sum_{u=1}^h
f_i(\mu_1,...,\mu_k)$$
 under $H$, where the intervals $[i_{u},i_{u-1}-1]$ for $u=1,...,h$ are given
in Remark 1.

Now focus on each $f_i$. Without loss of generality, focus on $f_1$,
then
$$f_1(\mu_1,...,\mu_k)=\sum_{j=i_1}^k n_j (\bar{Y}_j-\mu_j)^2=\sum_{j=i_1}^k n_i [(\bar{Y}_j-\hat{\mu}_j)^2
+2(\bar{Y}_j-\hat{\mu}_j)(\hat{\mu}_j-\mu_j)+(\hat{\mu}_j-\mu_j)^2].$$
Rearrange the terms above and note $\hat{\mu}_j=\bar{Y}_{i_1,k}$ for
$j\in[i_1,k]$, then
$$f_1(\mu_1,...,\mu_k)=\{\sum_{j=i_1}^k n_j[
(\bar{Y}_j-\hat{\mu}_j)^2+(\hat{\mu}_j-\mu_j)^2]\} +2\sum_{j=i_1}^k
n_j(\bar{Y}_j-\hat{\mu}_j)(-\mu_j)\stackrel{def}{=}I_1+I_2.$$ It is
obvious that $I_1$ is minimized at $\mu_j=\hat{\mu}_j$ for
$j\in[i_1,k]$; for $I_2$, apply Abel's partial summation formula and
obtain
$$I_2=2\sum_{j=k}^{i_1} n_j(\hat{\mu}_j-\bar{Y}_j)\mu_j=2
\sum_{j=k}^{i_1+1}d_j(\mu_j-\mu_{j-1})+d_{i_1}\mu_{i_1}=2
\sum_{j=k}^{i_1+1}d_j(\mu_j-\mu_{j-1}),$$ where $d_j=\sum_{v=k}^j
n_v (\hat{\mu}_v-\bar{Y}_v)\geq 0$ and $d_{i_1}=0$ due to the
definition of $i_1$. Also note $\mu_j\geq \mu_{j-1}$. Thus $I_2$ is
nonnegative and achieves its minimum at $\mu_{i_1}=...=\mu_k$.
Therefore, combining $I_1$ and $I_2$, we conclude
$f_1(\mu_1,...,\mu_k)$ is minimized at $\mu_j=\hat{\mu}_j$ for
$j\in[i_1,k]$.

Repeat the same argument on $f_2$ through $f_h$, each $f_u$ is
minimized at  $\mu_j=\hat{\mu}_j$ for $j\in[i_u, i_{u-1}-1]$.
Lastly, note $\hat{\mu}_j$ nondecreasing, we conclude
$\hat{\mu}_j^{mle}=\hat{\mu}_j$ for any $j\in[1,k]$.
 ~\raisebox{.5ex}{\fbox{}}

\bigskip

\begin{theorem}
\label{mle} For any $i\in[1,k]$,
\begin{equation}
\label{mle1}\hat{\mu_i}(\bar{Y}_1,...,\bar{Y}_k)=\hat{\mu}_i^{pava}(\bar{Y}_1,...,\bar{Y}_k).
\end{equation}
\end{theorem}

For a data set of $\{(\bar{Y}_i,n_i)\}_{i=1}^k$, the PAVA proceeds
as follows:

Step 0-PAVA). If $\bar{Y}_i$ is nondecreasing in $i$ for
$i\in[1,k]$, then $\hat{\mu}_i^{pava}=\bar{Y}_i$; otherwise, go to
the next step.

Step 1-PAVA). Pick any consecutive pair $(\bar{Y}_j, \bar{Y}_{j+1})$
with $\bar{Y}_j>\bar{Y}_{j+1}$, let $j_l$ be the smallest integer so
that $\bar{Y}_i=\bar{Y}_j$ for $i\in [j_l,j]$ and let $j_u$ be the
largest integer so that $\bar{Y}_i=\bar{Y}_{j+1}$ for $i\in
[j+1,j_u]$. For each $i\in[j_l,j_u]$, replace each $\bar{Y}_i$ by
$\bar{Y}_{j_l,j_u}(=\frac{\sum_{i=j_l}^{j_u}n_i\bar{Y}_i}{\sum_{i=j_l}^{j_u}n_i})$,
and then obtain a new data set of $\{(a_i,n_i)\}_{i=1}^k$, where
$a_i=\bar{Y}_i$ for $i\not \in [j_l,j_u]$ and
$a_i=\bar{Y}_{j_l,j_u}$ for $i \in [j_l,j_u]$. Note two facts:
$\bar{Y}_i$ is non-increasing for $i\in[j_l,j_u]$, and  the number
of different $a_i$'s is strictly less than that of $\bar{Y}_i$'s.

Step 2-PAVA) Repeat this process on $\{(a_i,n_i)\}_{i=1}^k$ until
all $a_i$'s are nondecreasing. Then $\hat{\mu}_i^{pava}=a_i$. Since
the number of different $a_i$'s is strictly less than that in the
previous step, this algorithm has to stop in a finite steps.
 ~\raisebox{.5ex}{\fbox{}}
\bigskip

\noindent {\bf Proof of Theorem~\ref{mle}.} When $\bar{Y}_i$ is
nondecreasing in $i\in[1,k]$, then (\ref{mle1}) is true due to Step
0-PAVA) and Remark 2. When $\bar{Y}_j>\bar{Y}_{j+1}$ for some $j$,
let $a_i$ and $[j_l,j_u]$ be given in Step 1-PAVA). It suffices to
show
\begin{equation}
\label{mleya}
\hat{\mu}_i(\bar{Y}_1,...,\bar{Y}_k)=\hat{\mu}_i(a_1,...,a_k),
\end{equation}
for any $i\in[1,k]$. Let $\cup_{u=1}^h [i_{u},i_{u-1}-1]$ be the
partition of $[1,k]$ given in Remark 1 using data
$\{(\bar{Y}_i,n_i)\}_{i=1}^k$. The integer $j_u$ has to belong to
one of these intervals in the partition, say $[i_{u_j},i_{u_j}-1]$.
Since $\bar{Y}_i$ is non-increasing on $[j_l,j_u]$ as shown in Step
1-PAVA), by Lemma~\ref{mon}, $i_{u_j}\leq j_l$. Thus $[j_l,j_u]$ is
a subset of $[i_{u_j},i_{u_j}-1]$, an interval in the partition.
 Let $\cup_{u=1}^{h'} [i_{u}',i_{u-1}'-1]$ be the
partition of $[1,k]$ given in Remark 1 but using data
$\{(a_i,n_i)\}_{i=1}^k$. Therefore, $[j_l,j_u]$ also has to be a
subset of one of these intervals.

Case I). If $[j_l,j_u]\subset [i_1,k]$, i.e., $u_j=1$, consider
$$a_{i,k}\stackrel{def}{=}\frac{\sum_{u=i}^k n_u a_u}{\sum_{u=i}^k
n_u}$$ for $i\in[1,k]$. Note $a_{i,k}=\bar{Y}_{i,k}$ for $i\not \in
[j_l+1,j_u]$ and $i_1'\leq j_l$, then $i_1'=i_1$. Therefore,
$$\hat{\mu}_i(\bar{Y}_1,...,\bar{Y}_k)=\frac{\sum_{i=i_1}^k
n_i\bar{Y}_i}{\sum_{i=i_1}^k n_i}=\frac{\sum_{i=i_1'}^k
n_ia_i}{\sum_{i=i_1'}^k n_i}=\hat{\mu}_i(a_1,...,a_k),$$ for any
$i\in[i_1,k]$. For any $i<i_1$, since $\hat{\mu}_i$ only depends on
$\bar{Y}_1$ through $\bar{Y}_{i_1-1}$(or $a_1$ through $a_{i_1-1}$)
and $\bar{Y}_i=a_i$, we conclude (\ref{mleya}).

Case II). If $[j_l,j_u]\subset [i_2,i_1-1]$, i.e., $u_j=2$, we only
need to show $i_1'=i_1$. Then, similar to Case I) above,
(\ref{mleya}) is established. To prove $i_1'=i_1$, first note
$$a_{i_1,k} =\bar{Y}_{i_1,k}\geq \bar{Y}_{u,k}= a_{u,k}$$ for any $u\geq i_1$.
So
\begin{equation}
\label{i1'} i_1'\leq i_1
\end{equation}
 by the definition of
$i_1'$.

Suppose $i_1'<i_1$, we will construct a contradiction. Note
$i_1'\not \in (j_l, j_u]$ because $a_i$ is non-increasing on $(j_l,
j_u]$. Thus, $$\bar{Y}_{i_1',k}=a_{i_1',k}\geq
a_{i_1,k}=\bar{Y}_{i_1,k},$$ a contradiction with the definition of
$i_1$. Hence $i_1'\geq i_1$. Together with (\ref{i1'}), we conclude
$i_1'=i_1$.

For the other cases of $u_j=3,...,h$, similar to Case II), we can
show $i_u=i_u'$ for all $u=2,...,h$. Hence, $h'=h$ and two
partitions, $\cup_{u=1}^h [i_{u},i_{u-1}-1]$ and $\cup_{u=1}^{h'}
[i_{u}',i_{u-1}'-1]$ are identical. Also note that $[j_l,j_u]$ is
contained in one interval $[i_{u_j},i_{u_j-1}-1]$, (\ref{mleya}) is
established. ~\raisebox{.5ex}{\fbox{}}
\bigskip

\noindent{\bf Remark 3}. Although $\hat{\mu}_i$ and
$\hat{\mu}_i^{pava}$ generated by two algorithms are identical,
there are several advantages of the SDMMSA over the PAVA. First, it
is clear from the SDMMSA that $\hat{\mu}_i$ is uniquely defined, but
not clear for $\hat{\mu}_i^{pava}$ from the PAVA, since the latter
needs to show $\hat{\mu}_i^{pava}$ must be the same no matter where
to start the algorithm, which is not obvious at all. Secondly,
$\hat{\mu}_j$ has a closed form, $\bar{Y}_{i_j, i_{j-1}-1}$, if
$j\in [i_{j_0}, i_{j_0-1}-1]$ for some $j_0$, where $i_{j_0}$ is
given in Remark 1, while $\hat{\mu}_i^{pava}$ does not. This fact is
important for establishing the stochastic ordering of $\hat{\mu}_i$
as shown in the next section. Thirdly, it was mentioned, for
example, in Robertson, Wright and Dykstra (1988, p.10) that
$\hat{\mu}_i^{pava}$ is the mle. To the best knowledge of the
authors, no rigorous proof has been given. With $\hat{\mu}_i$, we
proved $\hat{\mu}_i=\hat{\mu}_i^{mle}$ and
$\hat{\mu}_i=\hat{\mu}_i^{pava}$. Thus
$\hat{\mu}_i^{pava}=\hat{\mu}_i^{mle}$. Lastly, regarding the
computation, the SDMMSA is easier to code than the PAVA since at
each step of the SDMMSA a certain number of the final estimators
($\hat{\mu}_i$) are defined.
 ~\raisebox{.5ex}{\fbox{}}

\bigskip

\subsection{A stochastic ordering of $\hat{\mu}_i$.}
We provide another major result in this paper which establishes a
stochastic ordering for each $\hat{\mu}_i$ in terms of each of
$\mu_j$'s. Let
\begin{equation}
\label{muik} \hat{\mu}_{i,k}=\hat{\mu}_i(\bar{Y}_1,...,\bar{Y}_k)
\end{equation}
be the estimator of $\mu_i$ obtained from the sample $\{(\bar{Y}_i,
n_i)\}_{i=1}^k$ following Steps 1 and 2 in Section 2.1.
 So the distribution of $\hat{\mu}_{i,k}$ depends on $\mu_1$ through $\mu_k$ and
$\sigma$.
\begin{theorem}
\label{order} For each i and j in $[1,k]$,  $\hat{\mu}_{i,k}$, as a
function of $\bar{Y}_j$,  is nondecreasing when the other
$\bar{Y}_{j'}$'s are held fixed. Therefore, $\hat{\mu}_{i,k}$ is
stochastically nondecreasing in $\mu_j$ when the other $\mu_{j'}$'s
are held fixed. i.e., $P(\hat{\mu}_{i,k}>x)$ is a nondecreasing
function of $\mu_j$ for any real number $x$.
\end{theorem}
\noindent {\bf Proof of Theorem~\ref{order}.} We will prove the
monotonicity of $\hat{\mu}_{i,k}$ in each $\bar{Y}_j$ by induction
on $k$.

 First for
the case of $k=1$, $\hat{\mu}_{1,1}=\bar{Y}_1$ is nondecreasing in
$\bar{Y}_1$.

Assume that, for the case of $k=m$, $\hat{\mu}_{i,m}$ is
nondecreasing in $\bar{Y}_j$ for any $i$ and $j$ no larger than $m$.
Now consider the case of $k=m+1$. Following Step 1, $i_1$ depends on
$k$, so write it as $i_1(k)$, i.e., obtain $i_1(k)$ using
$(\bar{Y}_1, n_1)$ through  $(\bar{Y}_k, n_k)$. Claim
\begin{equation}
\label{i1m} i_1(m+1)=\left\{
\begin{array}{ll}
m+1, & \mbox{if} \,\ \bar{Y}_{m+1}>\hat{\mu}_{i_1(m),m},\\
i_1(m),   & \mbox{if} \,\ \bar{Y}_{m+1}\in
(\hat{\mu}_{i_2(m),m},\hat{\mu}_{i_1(m),m}],
\\
...\\
i_j(m),   & \mbox{if} \,\ \bar{Y}_{m+1}\in
(\hat{\mu}_{i_{j+1}(m),m},\hat{\mu}_{i_j(m),m}],\\
....\\
 i_h(m), &  \mbox{if} \,\ \bar{Y}_{m+1} \leq
\hat{\mu}_{i_h(m),m}.
\end{array}
           \right.
\end{equation}

When $\bar{Y}_{m+1}>\hat{\mu}_{i_1(m),m}(=\hat{\mu}_{m,m})$, for any
$j\in[1,m]$, note
$$ \bar{Y}_{j,m+1}=
\frac{n_{j,m} \bar{Y}_{j,m}+n_{m+1}\bar{Y}_{m+1}}{
n_{j,m}+n_{m+1}}\leq \frac{ n_{j,m}
\hat{\mu}_{m,m}+n_{m+1}\bar{Y}_{m+1}}{
n_{j,m}+n_{m+1}}<\bar{Y}_{m+1},$$ then $i_1(m+1)=m+1$.

When $\bar{Y}_{m+1}\in
(\hat{\mu}_{i_2(m),m},\hat{\mu}_{i_1(m),m}]=(\bar{Y}_{i_2(m),i_1(m)-1},\bar{Y}_{i_1(m),m}]$.
i) For $j\in[i_1(m),m+1]$, note $n_{j,m+1}-n_{i_1(m),m+1}\leq 0$,
then
\begin{eqnarray*}
&&\bar{Y}_{i_1(m),m+1}-\bar{Y}_{j,m+1}\\&=&\frac{n_{j,m+1}(n_{i_1(m),m}\bar{Y}_{i_1(m),m}+n_{m+1}\bar{Y}_{m+1})
-n_{i_1(m),m+1}(n_{j,m}\bar{Y}_{j,m}+n_{m+1}\bar{Y}_{m+1})}{n_{i_1(m),m+1}n_{j,m+1}}\\
&=&\frac{
n_{j,m+1}n_{i_1(m),m}\bar{Y}_{i_1(m),m}-n_{i_1(m),m+1}n_{j,m}
\bar{Y}_{j,m}+
(n_{j,m+1}-n_{i_1(m),m+1})n_{m+1}\bar{Y}_{m+1}}{n_{i_1(m),m+1}n_{j,m+1}}\\
&\geq& \frac{
n_{j,m+1}n_{i_1(m),m}\bar{Y}_{i_1(m),m}-n_{i_1(m),m+1}n_{j,m}
\bar{Y}_{i_1(m),m}+
(n_{j,m+1}-n_{i_1(m),m+1})n_{m+1}\bar{Y}_{i_1(m),m}}{n_{i_1(m),m+1}n_{j,m+1}}\\
&=&0,\end{eqnarray*} and conclude $i_1(m+1)\leq i_1(m).$

ii) For $j\in[1, i_1(m)-1]$, note
$$\bar{Y}_{j,m+1}=\frac{n_{j,i_1(m)-1}\bar{Y}_{j,i_1(m)-1}+n_{i_1(m),m}\bar{Y}_{i_1(m),m}+n_{m+1}\bar{Y}_{m+1}}{
n_{j,m+1}},$$ and $\bar{Y}_{j,i_1(m)-1}\leq
\bar{Y}_{i_2(m),i_1(m)-1}<\bar{Y}_{m+1}$, then
\begin{eqnarray*}
&&\bar{Y}_{i_1(m),m+1}-\bar{Y}_{j,m+1}\\
&=&
\frac{n_{j,i_1(m)-1}n_{i_1(m),m}\bar{Y}_{i_1(m),m}+n_{j,i_1(m)-1}n_{m+1}\bar{Y}_{m+1}-
n_{j,i_1(m)-1}(n_{i_1(m),m}+n_{m+1})\bar{Y}_{j,i_1(m)-1}}{n_{i_1(m),m+1}n_{j,m+1}}\\
&>&0,
\end{eqnarray*}
and conclude $i_1(m+1)> i_1(m)-1$. Therefore, combining i) and ii)
we obtain $i_1(m+1)=i_1(m)$ when $\bar{Y}_{m+1}\in
(\hat{\mu}_{i_2(m),m},\hat{\mu}_{i_1(m),m}]$.

For the other cases of $\bar{Y}_{m+1}$, ($\ref{i1m}$) can be
established in a similar way. Therefore, we conclude that
$\hat{\mu}_{i,m+1}$ depends on $\{(\bar{Y}_j,n_j)\}_{j=1}^{m+1}$
through $\{(\hat{\mu}_{j,m},n_j)\}_{j=1}^{m}$ and
$(\bar{Y}_{m+1},n_{m+1})$. So write
\begin{equation}
\label{mui}\hat{\mu}_{i,m+1}=\hat{\mu}_{i,m+1}(\hat{\mu}_{1,m},...,\hat{\mu}_{m,m},\bar{Y}_{m+1}).
\end{equation}
Also write $\hat{\mu}_{i,m+1}$ as
\begin{equation}
\label{muiyj}\hat{\mu}_{i,m+1}=\hat{\mu}_{i,m+1}(\bar{Y}_j),
\end{equation}
since the other $\bar{Y}_{j'}$ are fixed. We will use any one of the
above two notations whenever it is convenient.
 For $y<y'$, we need to show the monotonicity below
\begin{equation}
\label{thm1mono} \hat{\mu}_{i,m+1}(y) \leq \hat{\mu}_{i,m+1}(y'),
\end{equation}
which establishes the theorem, in the following two cases.

Case 1: $j=m+1$. Since $\hat{\mu}_{j',m}$ does not involve
$\bar{Y}_{m+1}$ for all $j'\leq m$, ($\ref{thm1mono}$) changes to
\begin{equation}
\label{thm1mono1}
\hat{\mu}_{i,m+1}(\hat{\mu}_{1,m},...,\hat{\mu}_{m,m},y) \leq
\hat{\mu}_{i,m+1}(\hat{\mu}_{1,m},...,\hat{\mu}_{m,m},y'),
\end{equation}
which is established in Lemma~\ref{lemmacase1} by noting
$\hat{\mu}_{i,m}$ is nondecreasing in $i\leq m$. Case 2: $j<m+1$.
Since $\hat{\mu}_{j',m}(y)\leq \hat{\mu}_{j',m}(y')$
 for all $j'\leq m$ by the induction assumption on the case of $k=m$, ($\ref{thm1mono}$) changes to
\begin{equation}
\label{thm1mono2}
\hat{\mu}_{i,m+1}(\hat{\mu}_{1,m}(y),...,\hat{\mu}_{m,m}(y),\bar{Y}_{m+1})
\leq
\hat{\mu}_{i,m+1}(\hat{\mu}_{1,m}(y'),...,\hat{\mu}_{m,m}(y'),\bar{Y}_{m+1}),
\end{equation}
which is established in Lemma~\ref{lemmacase2}. Therefore, the proof
of the monotonicity of $\hat{\mu}_{i,k}$ in each $\bar{Y}_j$ is
complete.

Since $\bar{Y}_j$'s are independent random variables, and each is
stochastically increasing in $\mu_j$, $\hat{\mu}_{i,k}$, as a
nondecreasing function of each $\bar{Y}_j$, is also stochastically
nondecreasing in $\mu_j$. See, for example, Alam and Rizvi (1966) or
Lemma 2 in Wu and Wang (2007). The proof of Theorem~\ref{order} is
complete.
 ~\raisebox{.5ex}{\fbox{}}

\bigskip

\begin{lemma}
\label{lemmacase1} Let
$\hat{\mu}_{i,m+1}(\bar{Y}_1,...,\bar{Y}_m,y)$ be the estimator
following Steps 1 and 2 on a date set
$\{(\bar{Y}_v,n_v)\}_{v=1}^{m+1}$ with $\bar{Y}_{m+1}=y$. Then
$$\hat{\mu}_{i,m+1}(\bar{Y}_1,...,\bar{Y}_m,y)\leq \hat{\mu}_{i,m+1}(\bar{Y}_1,...,\bar{Y}_m,y')$$
if $\bar{Y}_v$ is nondecreasing in $v\in[1,m]$ and $y<y'$.
\end{lemma}

\noindent{\bf Proof of Lemma~\ref{lemmacase1}.} Now write $i_1(m+1)$
introduced in ($\ref{i1m}$) as $i_1(y)$.

First  note $i_1(y)\leq i_1(y')$, which follows ($\ref{i1m}$) and
$y<y'$.

Secondly, claim
\begin{equation}
\label{i1y} i_1(y)=min\{v\in[1, m+1]: \bar{Y}_v\geq
\bar{Y}_{v,m+1}\}\stackrel{\mbox{denoted by}}{=}A.
\end{equation}
Note
\begin{equation}
\label{yv}
\bar{Y}_{v,m+1}-\bar{Y}_{v+1,m+1}=\frac{n_v(\bar{Y}_v-\bar{Y}_{v,m+1})}{n_{v+1,m+1}}=\frac{n_v(\bar{Y}_v-\bar{Y}_{v+1,m+1})}{n_{v,m+1}}.
\end{equation}
Therefore, $\bar{Y}_{v,m+1}$ is nonincreasing in $v$ when $v\geq A$.
Hence $i_1(y)\leq A$. It is obvious  that
$\bar{Y}_{A-1,m+1}<\bar{Y}_{A,m+1}$ following the first equality in
($\ref{yv}$). Since $\bar{Y}_v$ nondecreasing in $v\in [1,m]$,
$\bar{Y}_{v,m+1}$ is nondecreasing in $v$ when $v\leq A-1$. Thus
($\ref{i1y}$) is established.

Thirdly, a) when $i< i_1(y)$, both
$\hat{\mu}_{i,m+1}(\bar{Y}_1,...,\bar{Y}_m,y)$ and
$\hat{\mu}_{i,m+1}(\bar{Y}_1,...,\bar{Y}_m,y')$ are constructed
based on $\{\bar{Y}_v\}_{v=1}^{i_1(y')-1}$, a subset of
$\{\bar{Y}_v\}_{v=1}^m$ which is nondecreasing in $v$. Thus
$\hat{\mu}_{i,m+1}(\bar{Y}_1,...,\bar{Y}_m,y)=\bar{Y}_i=
\hat{\mu}_{i,m+1}(\bar{Y}_1,...,\bar{Y}_m,y')$ following Remark 2.

b) When $i\in [i_1(y), i_1(y')-1]$,
$\hat{\mu}_{i,m+1}(\bar{Y}_1,...,\bar{Y}_m,y)=\bar{Y}_{i_1(y),m+1}$
with $\bar{Y}_{m+1}=y$; while
$\hat{\mu}_{i,m+1}(\bar{Y}_1,...,\bar{Y}_m,y')=\bar{Y}_i\geq
\bar{Y}_{i_1(y),m+1}$ following ($\ref{i1y}$) and $\bar{Y}_v$
nondecreasing in $v\in[1,m]$.

c) When $i\in [i_1(y'), m+1]$,
$\hat{\mu}_{i,m+1}(\bar{Y}_1,...,\bar{Y}_m,y)=\bar{Y}_{i_1(y),m+1}$
with $\bar{Y}_{m+1}=y$; while
$\hat{\mu}_{i,m+1}(\bar{Y}_1,...,\bar{Y}_m,y')=\bar{Y}_{i_1(y'),m+1}$
with $\bar{Y}_{m+1}=y'$. Then
$\hat{\mu}_{i,m+1}(\bar{Y}_1,...,\bar{Y}_m,y)<
\hat{\mu}_{i,m+1}(\bar{Y}_1,...,\bar{Y}_m,y')$ due to $i_1(y)\leq
i_1(y')$, $\bar{Y}_v$ nondecreasing in $v\in[1,m]$, $y<y'$ and
($\ref{i1y}$). The proof is complete.
 ~\raisebox{.5ex}{\fbox{}}

\bigskip

\begin{lemma}
\label{lemmacase2} Let
$\hat{\mu}_{i,m+1}(\bar{Y}_1,...,\bar{Y}_m,\bar{Y}_{m+1})$ be the
estimator following Steps 1 and 2 on a date set
$\{(\bar{Y}_v,n_v)\}_{v=1}^{m+1}$. Then
$$\hat{\mu}_{i,m+1}(\bar{Y}_1,...,\bar{Y}_m,\bar{Y}_{m+1})\leq
\hat{\mu}_{i,m+1}(\bar{Y}_1',...,\bar{Y}_m',\bar{Y}_{m+1})$$ if
$\bar{Y}_v$ and $\bar{Y}_v'$ are both nondecreasing  and $\bar{Y}_v
\leq \bar{Y}_v'$ for $v\in[1,m]$.
\end{lemma}

\noindent{\bf Proof of Lemma~\ref{lemmacase2}.} Claim
\begin{equation}
\label{muv}
\hat{\mu}_{i,m+1}(\bar{Y}_1,...,\bar{Y}_{v-1},\bar{Y}_v,\bar{Y}_{v+1}..,\bar{Y}_m,\bar{Y}_{m+1})\leq
\hat{\mu}_{i,m+1}(\bar{Y}_1,...,\bar{Y}_{v-1},\bar{Y}_v',\bar{Y}_{v+1}..,\bar{Y}_m,\bar{Y}_{m+1})
\end{equation}
for any $v\in[1,m]$ if $\bar{Y}_1\leq...\leq \bar{Y}_{v-1}\leq
\bar{Y}_v\leq \bar{Y}_v'\leq \bar{Y}_{v+1}\leq ...\leq \bar{Y}_m$.
If ($\ref{muv}$) is true, then
$$\hat{\mu}_{i,m+1}(\bar{Y}_1,...,\bar{Y}_m,\bar{Y}_{m+1})\leq
\hat{\mu}_{i,m+1}(\bar{Y}_1,...,\bar{Y}_{m-1},
\bar{Y}_m',\bar{Y}_{m+1})\leq ...\leq
\hat{\mu}_{i,m+1}(\bar{Y}_1',...,\bar{Y}_m',\bar{Y}_{m+1}).$$

To show ($\ref{muv}$), now write $i_1(m+1)$ introduced in
($\ref{i1m}$) as $i_1(\bar{Y}_v)$(note $i_1(y)$ introduced in the
beginning of the proof of Lemma~\ref{lemmacase1} has a different
argument $y=\bar{Y}_{m+1}$). Following ($\ref{i1y}$) and the second
equality of ($\ref{yv}$), $i_1(\bar{Y}_v)\leq i_1(\bar{Y}_v')$.
Similar to the proof of Lemma~\ref{lemmacase1}, ($\ref{muv}$) can be
shown in three cases a) $i< i_1(y)$, b) $i\in [i_1(y), i_1(y')-1]$
and c) $i\in [i_1(y'), m+1]$, and the proof is complete.
 ~\raisebox{.5ex}{\fbox{}}

\bigskip

In short, in this section, we proposed the SDMMSA to generate
estimators for monotone normal means, $\mu_i$, proved that the
SDMMSA and the PAVA are equivalent, both generate the mle's, and the
distribution of the proposed estimator is stochastically
nondecreasing when $\mu_i$ goes larger. The last is to be used to
derive a test to detect the MED in the response-dose study as shown
in the next section.

\section{A step-up testing procedure to detect the MED.}

Now return to the problem of finding the minimum effective
dose(MED). First, we formulate this as a multiple test problem by
proposing a sequence of decreasing null hypotheses. Then a general
result that identifies the least favorable distribution is provided.
Finally, a sequence of increasing rejection regions of
level-$\alpha$ is constructed.

\subsection{Motivation.}

Let
\begin{equation}
\label{trans} X_i=\bar{Y}_i-\bar{Y}_0\,\ \mbox{and} \,\
\eta_i=\mu_i-\mu_0\,\ \mbox{for}\,\ i=1,...,k.
\end{equation}
Since the $MED$ is to be found, one should start the search from
$i=1$ instead of $i=k$.  Therefore, a step-up procedure seems more
reasonable than a step-down one. To establish $N=1$, some authors
(see, for example Hsu and Berger (1999)) compare $min\{X_j:j\geq
1\}$ with $\delta$ and claim $N=1$ if $min\{X_j:j\geq 1\}-\delta$ is
large in the unit of $S$. Roughly speaking, they use $min\{X_j:
j\geq i\}$ to estimate $\eta_i$. This does not fully utilize the
assumption of the monotonicity on means. So we propose using the
maximum likelihood estimator of $\eta_i$, denoted by
$\hat{\eta}_i\stackrel{def}{=}\hat{\mu}_i-\bar{Y}_0$, as a test
statistic, where $\hat{\mu}_i$ is given by the SDMMSA in Section
2.1. If $\hat{\eta}_1-\delta$ is larger than a multiple of $S$, then
claim $N=1$ and stop; otherwise compare $\hat{\eta}_2-\delta$ with
$S$. Repeat this process until we find an $N$ so that
$\hat{\eta}_{N}-\delta$ is much larger than $S$. If no such $N$ can
be found, then the $MED$ does not exist.

To identify $N$(MED), let
\begin{equation} \label{null}
{\cal C}=\{H_{0i}=\{\eta_i \leq \delta\}: i\in [1,k]\}
\end{equation}
be the set of null hypotheses of interest in this section. For each
$i\geq 1$, the alternative $H_{Ai}$ claims $\eta_i>\delta$. If a
certain $H_{Ai}$ is established, then $N\leq i$ due to the
monotonicity in $\mu_i$'s for $i\geq 1$. Therefore, $N$ should be
equal to the smallest $i$ so that $H_{Ai}$ is true. For the strong
control of the experimentwise error rate, it is clear that $H_{0i'}$
is a subset of $H_{0i}$ if $i < i'$ due to the monotonicity(i.e.,
$H_{0i}$ is decreasing). Therefore, ${\cal C}$ itself is closed
under the operation of intersection. The closed test procedure
(Marcus, Peritz and Gabriel, 1976) can be applied on ${\cal C}$ to
construct a step-up testing procedure with the experimentwise error
rate controlled at $\alpha$ in the strong sense (see, for example,
Hochberg and Tamhane (1987) for a definition) as long as a
level-$\alpha$ test is constructed for each $H_{0i}$. Let $R_i$ be a
rejection region  for $H_{0i}$ for any $i$ between 1 and k.  In
order to strongly control the experimentwise error rate, as well as
being powerful, region $R_i$ should satisfy the following two
properties:

$*)$    $R_i$ is of level $\alpha$, i.e.,
$\mbox{sup}_{\underline{\mu}\in
H_{0i}}P_{\underline{\mu}}(R_i)=\alpha.$

$**)$    $R_i$ is increasing in $i$. i.e. $R_i\subset R_{i'}$ if
$i<i'$. Thus $R_i=\cap_{\forall H_{0i'}\subset H_{0i}}
R_{i'}=\cap_{i'=i}^k R_{i'}$.

 \noindent If these two are satisfied, then the multiple
tests, which assert $H_{Ai}$ if and only if $R_i$ occurs, strongly
control the experimentwise error rate at level $\alpha$, which is
the main result of this section.

\bigskip

\subsection{A general result.}

\begin{theorem}
\label{thmg} Let $T(t_1,...,t_k)$ and $g_i(t_1,...,t_k)$ for
$i=1,...,k$ be  non-decreasing function for any  $t_i$ when the
other $t_j$'s are held constant. Also
\begin{equation}
\label{linear} g_i(ct_1+d,...,ct_k+d)=cg_i(t_1,...,t_k)+d
\end{equation}
for any constants $c>0$ and $d$.
 Then
\begin{equation}
\label{f}
f(\eta_1,...,\eta_k,\sigma)\stackrel{def}{=}ET(\frac{g_1(\bar{Y}_1,...,\bar{Y}_k)-\bar{Y}_0-\delta}{S},...,\frac{g_k(\bar{Y}_1,...,\bar{Y}_k)-\bar{Y}_0-\delta}{S})
\end{equation}
is nondecreasing in each $\eta_i$ when the other $\eta_j$ and
$\sigma$ are held constants.
\end{theorem}

\bigskip

\noindent{\bf Proof of Theorem~\ref{thmg}.} Due to ($\ref{linear}$),
we assume $\bar{Y}_0$ has a mean 0 and $\bar{Y}_i$ has a mean
$\eta_i(=\mu_i-\mu_0)$. Let $\phi(x)$ be the pdf of $N(0,1)$ and
$g_{\nu}(y)$ be the pdf of a $\chi^2$-distribution with
$\nu=\sum_{i=0}^k n_i-(k+1)$ degrees of freedom. Then
$$f(\eta_1,...,\eta_k,\sigma)=\int \int ET(\frac{g_1-x\frac{\sigma}{\sqrt{n_0}}-\delta}{
\sqrt{\frac{\sigma^2 y}{\nu}}},
...,\frac{g_k-x\frac{\sigma}{\sqrt{n_0}}-\delta}{
\sqrt{\frac{\sigma^2 y}{\nu}}})\phi(x)g_{\nu}(y)dx dy.$$ For each
fixed $x$ and $y$, let
$$T_{x,y}(\bar{Y}_1,...,\bar{Y}_k)=T(\frac{g_1-x\frac{\sigma}{\sqrt{n_0}}-\delta}{
\sqrt{\frac{\sigma^2 y}{\nu}}},
...,\frac{g_k-x\frac{\sigma}{\sqrt{n_0}}-\delta}{
\sqrt{\frac{\sigma^2 y}{\nu}}}),$$ which is non-decreasing in each
$\bar{Y}_i$ due to the monotonicity of $T$ and $g_i$'s. Therefore,
the conditional distribution of $T_{x,y}$ for given $x$ and $y$ is
stochastically nondecreasing in each $\eta_i$(see Lemma 2 in Wu and
Wang (2007)). Hence its conditional expectation
\begin{equation}
\label{inte}
ET_{x,y}=ET(\frac{g_1-x\frac{\sigma}{\sqrt{n_0}}-\delta}{
\sqrt{\frac{\sigma^2 y}{\nu}}},
...,\frac{g_k-x\frac{\sigma}{\sqrt{n_0}}-\delta}{
\sqrt{\frac{\sigma^2 y}{\nu}}}) \end{equation}
 is nondecreasing in each
$\eta_i$. So is $f$, the integral of ($\ref{inte}$).
 ~\raisebox{.5ex}{\fbox{}}
\bigskip

\noindent{\bf Remark 3}. Each
$g_i(\bar{Y}_1,...,\bar{Y}_k)\stackrel{def}{=}\hat{\mu}_i$ satisfies
($\ref{linear}$) and is nondecreasing in each $\bar{Y}_j$. We will
use this to construct step-up tests in the next section.
~\raisebox{.5ex}{\fbox{}}

\bigskip

\noindent{\bf Remark 4}. If define
$g_i(\bar{Y}_1,...,\bar{Y}_k)=\bar{Y}_i$ for $i\in[1,k]$ and a
sequence of
\begin{equation}
\label{hb} T^{HB}_j=I_{\{min_{\{i\in[j,k]\}}\{
(\bar{Y}_i-\bar{Y}_0-\delta)/\sqrt{1/n_i+1/n_0}\}>t_{\alpha,\nu}\}},
\end{equation}
for $j\in[1,k]$, then $g_i$ and $T^{HB}_j$ satisfy the conditions of
Theorem~\ref{thmg}. Hsu and Berger's step-down tests (1999) claim
$N$, the MED, to be $j_0$ if $T^{HB}_{j_0}=1$ but
$T^{HB}_{j_0-1}=0$. ~\raisebox{.5ex}{\fbox{}}
\bigskip

\subsection{The construction of step-up tests}

We first construct a rejection region $R^I_i$ with level $\alpha$
for each individual $H_{0i}$.
\begin{lemma}
\label{lem1} For a constant $c$, let
\begin{equation} \label{rIi}
R^I_{i,c}=\{\frac{\hat{\mu}_i-\bar{Y}_0-\delta}{S}>c\}.
\end{equation}
Then
\begin{equation} \label{typeI}
sup_{\underline{\bfb{\mu}} \in H_{0i}}
P_{\underline{\bfb{\mu}}}(R^I_{i,c})=
P_{\underline{\bfb{\mu}}_i}(R^I_{i,c}),
\end{equation}
 where $\underline{\bfb{\mu}}_i=(\mu_0,\mu_1,...,\mu_k)$ with
$\mu_1=...=\mu_i=\mu_0+\delta$ and $\mu_{i+1}=...=\mu_k=+\infty$.
Therefore, for any $\alpha\in (0,1)$, $R^I_{i,c}$, with
$c=c_{i,\alpha}$, defines a level-$\alpha$ test for $H_{0i}$,
 where $c_{i,\alpha}$ is the solution of
\begin{equation}\label{cut}
 P_{\underline{\bfb{\mu} }_i}(R^I_{i,c})=\alpha.
\end{equation}
\end{lemma}

\bigskip

\noindent{\bf Proof of Lemma~\ref{lem1}.} Let $T=I_{R^I_{i,c}}$.
Then Lemma~\ref{lem1} follows Theorem~\ref{thmg}.
~\raisebox{.5ex}{\fbox{}}

\bigskip

\noindent{\bf Remark 5}.
$c_{1,\alpha}=t_{\alpha,\nu}\sqrt{1/n_1+1/n_0}$ due to
$\hat{\mu}_1=\bar{Y}_1$ when
$\underline{\bfb{\mu}}=\underline{\bfb{\mu}}_1$.
 ~\raisebox{.5ex}{\fbox{}}

\bigskip

Region $R^I_{i,c}$ satisfies property *), but not property **) in
Section 2. To obtain more powerful multiple tests, we propose
\begin{theorem}\label{thm1}
For any integer $i\in [1,k]$ and for a sequence of nonnegative
constants $c_1$ through $c_i$,  let
\begin{equation}
\label{rc} R_{c_1,...,c_i}=\cup_{j=1}^i R^I_{j,c_j}=
\cup_{j=1}^i\{\frac{\hat{\mu}_{j}-\bar{Y}_0-\delta}{S}>c_j\}.
\end{equation}
 Then
\begin{equation}
\label{max}
 sup_{\underline{\bfb{\mu}}\in
H_{0i}}P_{\underline{\bfb{\mu}}}(R_{c_1,...,c_i})=P_{\underline{\bfb{\mu}}_i}(R_{c_1,...,c_i}).
\end{equation}
Therefore, for any $\alpha\in (0,1)$, $R_{c_1,...,c_i}$, with
$c_1=c_{1,\alpha}$ given in Remark 5 and $c_i$ determined
iteratively by solving

\begin{equation}
\label{cutm} P_{\underline{\bfb{\mu} }_i}(R_{c_1,...,,c_i})=\alpha,
\end{equation}
for $i=2,...,k$, defines a level-$\alpha$ test for $H_{0i}$.
\end{theorem}

\noindent{\bf Proof of Theorem~\ref{thm1}.} Let
$T=I_{R_{c_1,...,c_i}}$. Then Theorem~\ref{thm1} follows
Theorem~\ref{thmg}. ~\raisebox{.5ex}{\fbox{}}

\bigskip

\begin{theorem}
\label{thm3} Consider all hypotheses in ${\cal C}$ in $(\ref{null})$
with the following testing procedure:
\begin{equation}
\label{multi} \mbox{assert  $H_{Ai}$ (or not $H_{0i}$) if
$R_i\stackrel{def}{=}R_{c_1,...,c_i}$ occurs} \end{equation} for any
fixed $\alpha\in (0,1)$. Then the experimentwise error rate is at
most $\alpha$. i.e., the probability of making at least one
incorrect assertion is at most $\alpha$.
\end{theorem}

\noindent{\bf Proof of Theorem~\ref{thm3}}. The proof is trivial if
one notices that $H_{0i}$ is decreasing in $i$, and $R_i$ is of
level-$\alpha$ and is increasing in $i$. Then Theorem~\ref{thm3}
follows the closed test procedure by Marcus, Peritz and Gabriel
(1976).
 ~\raisebox{.5ex}{\fbox{}}

\bigskip

\noindent{\bf Remark 6.} Region $R_i$ is increasing in $i$. Then
$R_i$ satisfies properties *) and **).
 ~\raisebox{.5ex}{\fbox{}}

\bigskip

\noindent{\bf Remark 7.} When the design is balanced, region $R_1$
contains the set of $\{T^{HB}_1=1\}$, on which Hsu and Berger's test
(1999) claims the MED=1. Therefore, the proposed test is uniformly
more powerful than Hsu and Berger's when the MED=1.
 ~\raisebox{.5ex}{\fbox{}}

\bigskip

\noindent{\bf Example 1(continued).}
 The sample standard deviation
$S=7.751$, $t_{0.05, 50}=1.676$. We compare the new step-up
procedure with the step-up procedure SU1P, the step-down Williams
procedure and step-down procedure SD1P in Tamhane et al (1996) and
the DR method in Hsu and Berger (1999). For illustration, $\delta
=6.5$. From Table 1 in Dunnett and Tamhane (1992) we have the
critical values for the step-up procedure SU1P $c_1=1.645,
c_2=1.933, c_3=2.071, c_4=2.165, c_5=2.237, c_6=2.294, c_7=2.342,
c_8=2.382$ (we treat $df=50$ as $df = \infty$). The SU1P procedure
infers $\widehat{MED} = 5$.
 The Williams procedure
has the $\bar{t}$ statistics: $\bar{t}_1=-1.810, \bar{t}_2=-0.961,
\bar{t}_3=0.313, \bar{t}_4=1.899, \bar{t}_5=5.788,
\bar{t}_6=\bar{t}_7=\bar{t}_8=9.334, \bar{t}_9=9.877$. The Williams
statistics are compared with the following critical values (taken
from Williams (1971)) in a step-down manner: $c_1=1.675, c_2=1.755,
c_3=1.780, c_4=1.790, c_5=1.795, c_6=1.800, c_7=1.805, c_8=1.805,
c_9=1.810.$. The Williams procedure infers $\widehat{MED} = 4$. By
simulation with 1,0000 repetition, the critical values of the new
statistic are $c_1=0.968, c_2=1.022, c_3=1.046, c_4=1.046,
c_5=1.034, c_6=1.043, c_7=1.044, c_8=1.047, c_9=1.030$,
respectively. Also $(\hat{\mu}_i-\bar{Y}_0-\delta)/S$ for
$i=1,2,3,4$ given in (\ref{rc}) are  -1.045,
     -0.555,
      0.181, and
       1.097, respectively. Thus the new step-up procedure concludes  $\widehat{MED} = 4$. So does
        Hsu and Berger (1999)'s DR method.

\section{Discussion.}

In this paper, we propose an alternative, SDMMSA, for the widely
used PAVA. Although the two are equivalent, the SDMMSA is important
by itself since it is easily coded and is well defined. Also a
stochastic ordering of the estimators for the monotone normal means
is established based on the SDMMSA. As one of its applications, a
step-up test procedure is proposed to identify the MED. It strongly
controls the experimentwise error rate, and is powerful to detect
the MED, especially when the true MED is small.

\section*{Acknowledgments.} We thank Professor Roger Berger for his helpful comments.

\newpage

\begin{center}

{TABLE 1. Sample dose-response data in Example 1}\\

\begin{tabular}{ccccccc}\hline \hline

 Dosage&Sample &$\bar{Y}_i$&SD&Index&$\hat{\mu}_i$& $i_{j}$\\

 (mg/kg)&size& &response&&&\\\hline

 0&6&25.5&2.6&0&-&-\\

 0.5&6&23.9&4.0&1&23.9&$i_7=1$\\

 1.0&6&27.7&3.3&2&27.7&$i_6=2$\\

 1.5&6&33.4&2.3&3&33.4&$i_5=3$\\

 2.0&6&40.5&10.5&4&40.5&$i_4=4$\\

 2.5&6&57.9&9.9&5&57.9&$i_3=5$\\

 3.0&6&74.4&14.6&6&73.77&$i_2=6$\\

 3.5&6&73.4&7.6&7&73.77&-\\

 4.0&6&73.5&4.5&8&73.77&-\\

 4.5&6&76.2&7.9&9&76.2&$i_1=9$\\\hline \hline

\end{tabular}

\end{center}
\end{document}